\documentclass[a4paper]{amsart}

\usepackage[english]{babel}
\usepackage[utf8x]{inputenc}
\usepackage[T1]{fontenc}
\usepackage{subdepth}
\usepackage{cite}

\usepackage{amsthm,amsmath,amssymb}
\usepackage{tikz-cd}
\usepackage{amsfonts}
\usepackage{graphicx}
\usepackage[colorlinks=true, allcolors=blue]{hyperref}
\usepackage{mathdots}
\usepackage{csquotes}
\usepackage{xfrac}

\theoremstyle{plain}
\newtheorem{thm}{Theorem}[section]
\newtheorem{lem}[thm]{Lemma}
\newtheorem{prop}[thm]{Proposition}
\newtheorem{cor}[thm]{Corollary}
\newtheorem*{thm*}{Theorem}
\numberwithin{equation}{section}
\newtheorem{thmx}{Theorem}

\theoremstyle{definition}
\newtheorem{example}[thm]{Example}
\newtheorem{defn}[thm]{Definition}

\theoremstyle{remark}

\newtheorem{rmk}[thm]{Remark}

\newcommand*{\suchthat}{\;\ifnum\currentgrouptype=16 \middle\fi|\;}
\newcommand{\Z}{\mathbb{Z}}

\newcommand{\acts}{\mathbin{\curvearrowright}}
\newcommand{\rightacts}{\mathbin{\curvearrowleft}}

\renewcommand{\H}{\mathrm{H}}
\newcommand{\K}{\mathrm{K}}
\DeclareMathOperator{\Ab}{\mathbf{Ab}}
\newcommand{\lMod}[1]{#1\text{-}\mathbf{Mod}} 
\newcommand{\rMod}[1]{\mathbf{Mod}\text{-}#1}

\newcommand{\cat}[1]{\textsf{#1}}
\newcommand{\bb}[1]{\mathbb{#1}}

\DeclareMathOperator{\Tor}{Tor}
\DeclareMathOperator{\End}{End}

\newcommand{\cal}[1]{\mathcal{#1}}
\newcommand{\bcdot}{\boldsymbol{\cdot}}

\DeclareMathOperator{\id}{id}

\DeclareMathOperator{\algspan}{span}
\DeclareMathOperator{\ev}{ev}

\DeclareMathOperator{\Ind}{Ind}
\DeclareMathOperator{\Res}{Res}

\DeclareMathOperator{\Coinv}{Coinv}

\newcommand{\defeq}{:=}
\newcommand{\orho}{\overline{\rho}}
\newcommand{\cs}{\mathrm{C}^\ast}
\DeclareMathSymbol{\mhyphen}{\mathord}{AMSa}{"39}

\makeatletter
\renewenvironment{cases}[1][l]{\matrix@check\cases\env@cases{#1}}{\endarray\right.}
\def\env@cases#1{%
  \let\@ifnextchar\new@ifnextchar
  \left\lbrace\def\arraystretch{1.2}%
  \array{@{}#1@{\quad}l@{}}}
\makeatother

\newcommand\extrafootertext[1]{%
    \bgroup
    \renewcommand\thefootnote{\fnsymbol{footnote}}%
    \renewcommand\thempfootnote{\fnsymbol{mpfootnote}}%
    \footnotetext[0]{#1}%
    \egroup
}

\title{Ample groupoid homology and \'etale correspondences}
\author{Alistair Miller}

\begin{document}

\maketitle
\begin{abstract}
We show that \'etale correspondences between ample groupoids induce homomorphisms of homology groups. To complement this we explore the module categories of ample groupoids. We construct an induction-restriction adjunction for subgroupoids, which generates a procedure for building resolutions of arbitrary groupoid modules. These resolutions can be used to work with the Tor picture of groupoid homology, enabling explicit descriptions of the maps in homology induced by \'etale correspondences.
\end{abstract}

\extrafootertext{This work contains part of the author's PhD thesis, which was supported by the Engineering and Physical Sciences Research Council (EPSRC) through a doctoral studentship. The author has also been supported by the European Research Council (ERC) under the European Union's Horizon 2020 research and innovation programme (grant agreement No. 817597) and by the Independent Research Fund Denmark through the Grant 1054-00094B.}

\section{Introduction}

This paper is focused on a homology theory for ample groupoids. Introduced in the \'etale setting \cite{CraMoe00}, the theory for ample groupoids in particular garnered attention from those in $\cs$-algebras and in topological dynamics after a series of papers by Matui \cite{Matui12, Matui15, Matui16}. Here he advanced two conjectures for an ample groupoid $G$, the HK conjecture relating the homology $\H_*(G)$ to the K-theory $\K_*(\cs_r(G))$ of the groupoid $\cs$-algebra, and the AH conjecture relating $\H_*(G)$ to the topological full group of $G$. Although there are counterexamples to the HK conjecture \cite{Scarparo20, Deeley21}, both conjectures have been verified for large classes of ample groupoids (see \cite{FKPS19, BDGW23} and \cite{Scarparo20, Ortega20, NylOrt21a, NylOrt21b, Scarparo22, Li22}) and have led to the discovery of deeper connections of groupoid homology to operator K-theory \cite{ProYam22} and topological full groups \cite{Li22}.

Many computations in groupoid homology utilise its invariance under Morita equivalences \cite{CraMoe00}. We work more generally with \'etale correspondences, a notion of morphism for groupoids which capture not only Morita equivalences, but also \'etale homomorphisms and actors/algebraic morphisms. There are other examples, such as a correspondence $\Omega_S \colon S \ltimes E^\times \to G_S$ which we associate to an inverse semigroup $S$. \'Etale correspondences were introduced to serve as models for $\cs$-correspondences of the associated $\cs$-algebras \cite{Holkar14, Holkar17b, Holkar17a, Albandik15, AKM22}. $\cs$-correspondences are the basis for Kasparov cycles, and a proper \'etale correspondence of \'etale groupoids induces a map in the K-theory of the associated $\cs$-algebras. In accordance with the similarities between groupoid homology and K-theory, we develop the functoriality of groupoid homology with respect to \'etale correspondences. This generalises the Morita invariance of groupoid homology for ample groupoids.

\begin{thmx}[Corollary \ref{main corollary}]\label{short main theorem}
Let $G$ and $H$ be ample groupoids. Then a proper \'etale correspondence $\Omega \colon G \to H$ induces homomorphisms
\[ \H_*(\Omega) \colon \H_*(G) \to \H_*(H) \]
of homology groups.
\end{thmx}

The map $\H_*(\Omega) \colon \H_*(G) \to \H_*(H)$ is designed to be able to describe homomorphisms of homology groups which may appear in other contexts. For example, in \cite{ProYam22} the homology $\H_*(G)$ is recovered from equivariant KK-theoretic data in order to connect it to the K-theory $\K_*(C^*_r(G))$. The K-theoretic functoriality of \'etale correspondences is extended to the level of equivariant KK-theory in \cite{Miller23a}. To show that $\H_*(\Omega) \colon \H_*(G) \to \H_*(H)$ is compatible with this KK-theoretic functoriality via the recovery process, the description of $\H_*(\Omega)$ must be flexible enough to apply in this situation. This means working with groupoid modules which appear as K-theory groups of $\cs$-algebras equipped with the action of a groupoid.

Working with the category $\lMod G$ of $G$-modules we have a flexible description of the homology groups $\H_*(G)$ as $\Tor^G_*(\Z[G^0],\Z[G^0])$ \cite{BDGW23, Li22}. Concretely this means that any projective or even left $G$-acyclic\footnote{A $G$-module $P$ is left $G$-acyclic if $\H_i(G;P) = 0$ for each $i \geq 1$. This property is weaker than both projectivity and flatness.} resolution $P_\bullet \to \Z[G^0]$ of $G$-modules computes $\H_*(G)$ as $\H_*( \Z[G^0] \otimes_G P_\bullet)$. For example, using the bar resolution $\Z[G^{\bullet +1}] \to \Z[G^0]$ recovers Matui's definition of $\H_*(G)$ \cite[Definition 3.1]{Matui12}. The value of Corollary \ref{main corollary} is not only the existence of $\H_*(\Omega)$ but also its uniqueness: it allows us to recognise when a map from $\H_*(G)$ to $\H_*(H)$ is equal to $\H_*(\Omega)$, no matter which resolutions are used to model $\H_*(G)$ and $\H_*(H)$.

In the setting of an inverse semigroup $S$, the correspondence $\Omega_S \colon S \ltimes E^\times \to G_S$ induces an isomorphism in homology (Example \ref{inverse semigroup homology isomorphism}). This computes the homology $\H_*(G_S)$ of the universal groupoid $G_S$ as 
\[ \H_*(G_S) \cong \bigoplus_{[e] \in S \backslash E^\times} \H_*(S_e), \]
where $S_e = \{ s \in S \suchthat s^* s = e = s s^* \}$ is the stabiliser subgroup of $S$ at a non-zero idempotent $e \in E^\times$.

The module-theoretic approach to Theorem \ref{short main theorem} extends naturally to the homology $\H_*(G;M) = \Tor^G_*(\Z[G^0],M)$ of an ample groupoid $G$ with coefficients in a $G$-module $M$. The properness of the \'etale correspondence $\Omega \colon G \to H$ is used only to define a $G$-module map from $\Z[G^0]$ to the $G$-module $\Ind_\Omega \Z[H^0]$ induced by $\Omega$. In general, an \'etale correspondence $\Omega \colon G \to H$ induces maps in homology with coefficients as follows.
\begin{thmx}[Theorem \ref{inducedmapinhomology}]\label{main theorem}
Let $G$ and $H$ be ample groupoids, let $M$ be a $G$-module and $N$ an $H$-module. Then for each \'etale correspondence $\Omega \colon G \to H$ and $G$-module map $f \colon M \to \Ind_\Omega N$, there are induced homomorphisms
\[ \H_*(\Omega;f) \colon \H_*(G;M) \to \H_*(H;N) \]
of homology groups.
\end{thmx}

To complement our approach which enables the use of arbitrary left $G$-acyclic resolutions, we construct a canonical left $G$-acyclic resolution of any $G$-module. This derives from an induction-restriction adjunction, reflecting the situation for group modules.
\begin{thmx}[Proposition \ref{induction restriction}]
Let $G$ be an ample groupoid and let $H \subseteq G$ be an open subgroupoid. Then there is a restriction functor $\Res^H_G \colon \lMod G \to \lMod H$ and an induction functor $\Ind^H_G \colon \lMod H \to \lMod G$ which is left adjoint to $\Res^H_G$.
\end{thmx}
Taking the unit space $X = G^0$, the adjunction $\Ind^G_X \dashv \Res^X_G$ generates in a standard way the bar resolution
\[ \cdots \to \Z[G^{n+1}] \otimes_X M \to \cdots \to \Z[G] \otimes_X M \to M  \]
of any $G$-module $M$. This is left $G$-acyclic by a groupoid version of Shapiro's Lemma (Lemma \ref{Shapiros lemma}), and it recovers the chain complex given in \cite[Section 1.3]{ProYam22} to define $\H_*(G;M)$ after taking the coinvariants $\Z[G^0] \otimes_G -$.

The author would like to thank Jamie Gabe, Xin Li and Owen Tanner for helpful comments and discussions, and Ben Steinberg for pointing out an error in a previous version of the manuscript.

\section{Modules over an ample groupoid}

An \textit{ample} groupoid is an \'etale groupoid whose unit space is locally compact, Hausdorff and totally disconnected. The space of arrows is then automatically totally disconnected and locally LCH, but not necessarily Hausdorff. Following for example \cite{BDGW23} and \cite{Li22}, we take a module approach to the homology of an ample groupoid. For a totally disconnected locally LCH space $X$ and an open Hausdorff subspace $U \subseteq X$, we include $C_c(U, \Z)$ into the abelian group of integer valued functions on $X$ by setting an element of $C_c(U, \Z)$ to be $0$ outside $U$. We write $\Z [X]$ for the span of $C_c(U,\Z)$ across all open Hausdorff subspaces $U \subseteq X$. For any open Hausdorff cover $\mathcal U$ of $X$, the span of $C_c(U, \Z)$ over $U \in \mathcal U$ is equal to $\Z[X]$.

\begin{defn}[Groupoid ring]
For an ample groupoid $G$, the \textit{groupoid ring} is the abelian group $\Z[G]$ with a multiplication given by convolution: for $\xi, \eta \in \Z[G]$ and $g \in G$ we have 
\[ \xi * \eta (g) = \sum_{h \in G_{r(g)}} \xi(h^{-1}) \eta(hg). \]
\end{defn}

In this document, a ring need not be commutative nor even unital. In place of unitality, $\Z[G]$ is \textit{locally unital}. This means that for any finite collection $\xi_1, \dots, \xi_n$ of elements in $\Z[G]$, there is an idempotent $e \in \Z[G]$ such that $e\xi_i = \xi_i e = \xi_i$ for each $i$. In this case the idempotent may be taken to be the indicator function $\chi_U$ on a compact open set $U \subseteq G^0$. For a locally unital ring $R$, we require our (left) $R$-modules $M$ to be \textit{non-degenerate} or \textit{unitary} in the sense that $R M = M$. The categories of left and right $R$-modules are written $\lMod R$ and $\rMod R$ respectively. When $R = \Z[G]$, we refer to $R$-modules as $G$-modules and write $\lMod G$ and $\rMod G$ for the left and right module categories. 

For continuous maps $\sigma \colon Y \to X$ and $\rho \colon Z \to X$ we write $Y \times_{\sigma, X, \rho} Z$ for the \textit{fibre product} $\{ (y,z) \in Y \times Z \suchthat \sigma(y) = \rho(z)\}$ of $Y$ and $Z$ over $X$. If $\sigma$ and $\rho$ are understood we may write $Y \times_X Z$.

\begin{defn}[$G$-space]
Let $G$ be an \'etale groupoid. A (left) \textit{$G$-space} $X$ consists of a topological space equipped with:
\begin{itemize}
\item a continuous map $\tau \colon X \to G^0$ called the \textit{anchor map},
\item a continuous map $\alpha \colon G \times_{s,G^0,\tau} X \to X$ called the \textit{action map}. We usually write $g \bcdot x$ for $\alpha(g,x)$.
\end{itemize}
This describes a $G$-space if whenever $g, h \in G$ and $x \in X$ satisfy $s(g) = r(h)$ and $s(h) = \tau(x)$, we have $\tau(h \bcdot x) = r(h)$ and $gh \bcdot x = g \bcdot (h \bcdot x)$. 
\end{defn}

\begin{example}[$G$-space module]\label{G-space module}
Let $G$ be an ample groupoid and let $X$ be a totally disconnected locally LCH (left) $G$-space with anchor map $\tau \colon X \to G^0$. Then the abelian group $\Z[X]$ is a $G$-module with $\Z[G]$-action given by
\begin{align*}
\Z[G] \times \Z[X] & \to \Z[X] \\
(\xi, m) & \mapsto \xi \bcdot m \\
& \quad \; \; x \mapsto \sum_{g \in G_{\tau(x)}} \xi(g^{-1}) m(g \bcdot x).
\end{align*} 
Similarly, for a totally disconnected locally LCH right $G$-space $Z$ with anchor map $\pi \colon Z \to G^0$, the abelian group $\Z[Z]$ is a right $G$-module via the map
\begin{align*}
\Z[Z] \times \Z[G] & \to \Z[Z] \\
(m,\xi) & \mapsto m \bcdot \xi \\
& \quad \; \; z \mapsto \sum_{g \in G^{\pi(z)}} m(z \bcdot g) \xi(g^{-1}).
\end{align*}
Given a $G$-equivariant local homeomorphism $f \colon X \to Y$ of (left) $G$-spaces, we functorially obtain a $G$-module homomorphism $f_* \colon \Z[X] \to \Z[Y]$, which at $\xi \in \Z[X]$ and $y \in Y$ is given by
\[ f_* (\xi) (y) = \sum_{x \in f^{-1}(y)} \xi(x). \]
The canonical left and right actions of $G$ on its unit space $G^0$ give rise to modules with abelian group $\Z[G^0]$ which we call the \textit{trivial} left and right $G$-modules. 
\end{example}
Note that we do not require our $G$-spaces to have an \'etale anchor map as in \cite{BDGW23}. To work with the abelian group $\Z[X]$ of a totally disconnected locally LCH space $X$ we make use of the following rephrasing of \cite[Lemma 2.2]{Li22}, which can for example be used to construct actions of $\Z[G]$.

\begin{lem}\label{basiscompatibility}
Let $X$ be a totally disconnected locally LCH space and let $A$ be an abelian group. Suppose that $\cal O$ is an open Hausdorff cover of $X$ and let $\cal U$ be the set of compact open subsets of elements of $\cal O$. Let $\varphi \colon \cal U \to A$ be a function such that $\varphi(U_1) + \varphi(U_2) = \varphi(U)$ whenever $U_1 \sqcup U_2 = U$ with $U_1$, $U_2$ and $U$ in $\cal U$.

Then $\varphi \colon \cal U \to A$ extends uniquely to a group homomorphism $\hat{\varphi} \colon \Z[X] \to A$ such that $\hat{\varphi}(\chi_U) = \varphi(U)$ for each $U \in \cal U$.
\end{lem}

The categories $\lMod R$ and $\rMod R$ of left and right modules over a locally unital ring $R$ retain much of the behaviour of the unital setting. We may take the tensor product $M \otimes_R N$ of a right $R$-module $M$ and a left $R$-module $N$, and a (left) $R$-module is \textit{flat} if the functor $- \otimes_R M \colon \rMod R \to \Ab$ preserves exact sequences, or equivalently preserves injectivity. An $R$-module $P$ is \textit{projective} if every morphism out of $P$ lifts through surjective morphisms. Direct sums of copies of $R$ need no longer be freely generated by a set of elements, and $R$ need not be projective as a module over itself\footnote{In \cite[Section 2.3]{BDGW23} this is pointed out for $\Z[X]$ where $X$ is the reals equipped with a rational sequence topology.}, although it is always flat. Instead, for each idempotent $e \in R$, the module $Re$ is projective \cite[Remark 2.9]{BDGW23}. It is then straightforward to see that $\lMod R$ has enough projectives, as for any $R$-module $M$, there is a direct sum $P$ of modules of the form $Re$ and a surjective morphism $P \to M$.

\begin{rmk}\label{multiplier ring}
We may often carry over results from the unital setting by using the multiplier ring $M(R) \defeq \End_{\rMod R}(R)$ of a locally unital ring $R$. This is a unital ring containing $R$, and $R$-modules may be identified with $M(R)$-modules $A$ such that $R \bcdot A = A$. Under this identification, an $R$-module is projective/flat if and only if it is projective/flat as an $M(R)$-module.
\end{rmk}

\begin{defn}
Let $G$ be an \'etale groupoid and let $X$ be a topological space. An action $G \acts X$ is \textit{basic} if the map 
\begin{align*}
 G \times_{G^0} X & \to X \times_{G \backslash X} X  \\
 (g,x) & \mapsto (g \bcdot x, x) 
\end{align*}
is a homeomorphism.
\end{defn}

A basic action $G \acts X$ is free, and is proper if and only if the orbit space $G \backslash X$ is Hausdorff \cite[Proposition 2.19]{AKM22}. Conversely, any free proper action is basic. The orbit map $q \colon X \mapsto G \backslash X$ of a basic action is \'etale \cite[Lemma 2.12]{AKM22}, so if $X$ is locally LCH then so is $G \backslash X$. For a basic $G$-space $X$ which is further \textit{\'etale} in the sense that its anchor map is \'etale, a \textit{slice} $U \subseteq X$ is an open Hausdorff subspace on which the quotient map and the anchor map are injective. This defines a partial homeomorphism from $G^0$ to $G \backslash X$. Slices for the left multiplication action $G \acts G$ are exactly the open bisections.

For the next proposition we will utilise Steinberg's picture of $G$-modules for an ample groupoid $G$ in terms of $G$-sheaves \cite{Steinberg14}. For each $G$-module $M$ there is a $G$-sheaf $\cal M$ along with a natural isomorphism $M \cong \Gamma_c(G^0, \cal M)$ of $G$-modules. The fibre $M_x$ of $\cal M$ at $x \in G^0$ is constructed as a direct limit over the compact open subsets $U, V \subseteq G^0$ containing $x$. Given $x \in V \subseteq U$, multiplication by $\chi_V$ defines a homomorphism $\chi_U \bcdot M \to \chi_V \bcdot M$, and we have
\[ M_x =  \varinjlim_{x \in U} \chi_U \bcdot M. \]
As $M$ is the union of $\chi_U \bcdot M$ for such $U$, we may assign to each $m \in M$ its equivalence class $m_x \in M_x$. The topology on $\cal M$ ensures that for each $m \in M$ the section 
\begin{align*}
G^0 & \to \cal M \\
x & \mapsto m_x 
\end{align*}
is continuous and compactly supported. If $m_x = n_x$ for some $m, n \in M$ and $x \in G^0$ then by construction there is a compact open neighbourhood of $x$ on which $m$ agrees with $n$. A $G$-map $f \colon M \to N$ induces a morphism of the associated $G$-sheaves $(f_x)_{x \in G^0} \colon \cal M \to \cal N$, and $f$ is injective if and only if $f_x \colon M_x \to N_x$ is injective for each $x \in G^0$. 
\begin{rmk}
Any sheaf $\cal M$ over a totally disconnected LCH space $X$ is automatically c-soft in that each continuous section $f \colon K \to \cal M$ on a compact subspace $K \subseteq X$ extends to a continuous section $\tilde f \colon X \to \cal M$. 
\end{rmk}

For any ample groupoid $G$ and $n \geq 1$ the $G$-space module $\Z[G^n]$ is flat \cite[Proposition 2.4]{Li22}. This result may be adapted to any basic \'etale $G$-space:
\begin{prop}\label{flat module result}
Let $G$ be an ample groupoid and let $X$ be a basic \'etale right $G$-space. Then $\Z[X]$ is a flat $G$-module.

\begin{proof}
Let $\sigma \colon X \to G^0$ be the anchor map and $q \colon X \to X/G$ the orbit map. Let $f \colon A \to B$ be an injective $G$-module map. This implies that the fibre $f_y \colon A_y \to B_y$ is injective for each $y \in G^0$. We need to check that $\id \otimes f \colon \Z[X] \otimes_G A \to \Z[X] \otimes_G B$ is injective, so let $a \in \ker \id \otimes f$. Then there are compact slices $U_i \subseteq X$ and elements $a_i \in A$ indexed by a finite set $I$ such that
\[ a = \sum_{i \in I} \chi_{U_i} \otimes a_i.\] 
Let $\{O_1,\dots,O_n\}$ be a finite set of compact open Hausdorff subsets of $X/G$ such that for each $i \in I$, there is $1 \leq k \leq n$ such that $q(U_i) \subseteq O_k$. We proceed to show that $a = 0$ by induction on $n$. 

If $n = 1$, then $q(U_i) \subseteq O_1$ for each $i \in I$. We may assume that the $q(U_i)$ are pairwise disjoint because of disjointification: for each non-empty $J \subseteq I$ let 
\[V_J = \bigcap_{i \in J} q(U_i) \cap \bigcap_{i \in I \setminus J} (O_1 \setminus q(U_i)),\]
and for $i \in J$ set $U'_{i,J} = q^{-1}(V_J) \cap U_i$. Then the $V_J$ are pairwise disjoint and for each $i \in I$, we have 
\[U_i = \bigsqcup_{i \in J \subseteq I} U'_{i,J}. \]
For each $J$ pick an element $i_J \in J$. For each $i \in J$ we have $q(U'_{i,J}) = V_J = q(U'_{i_J,J})$ so because the action $X \rightacts G$ is basic there is a compact open bisection $W_{i,J} \subseteq G$ such that $U'_{i,J} = U'_{i_J,J} \bcdot W_{i,J}$. Setting 
\[b_J = \sum_{i \in J } \chi_{W_{i,J}} \bcdot a_i, \]
then by construction we have
\[
a = \sum_{i \in I} \chi_{U_i} \otimes a_i  
= \sum_{\emptyset \ne J \subseteq I} \sum_{i \in J} \chi_{U'_{i,J}} \otimes a_i    
= \sum_{\emptyset \ne J \subseteq I} \chi_{U'_{i_J, J}} \otimes b_J.
\]
Thus we may assume that the $q(U_i)$ are pairwise disjoint for $i \in I$. For each $x \in X$ there is an evaluation map for $A$ given by
\begin{align*}
\ev_x \colon \Z[X] \otimes_G A & \to A_{\sigma(x)} \\
\xi \otimes a & \mapsto \sum_{g \in G^{\sigma(x)}} \xi(x \bcdot g) (g^{-1} \bcdot a_{r(g)})
\end{align*}
and similarly for $B$. These are compatible with $f \colon A \to B$ in that $\ev_x \circ (\id \otimes f) = f_{\sigma(x)} \circ \ev_x$. For $x \in U_i$, pairwise disjointness of the orbits of the slices implies that $\ev_x(a) = (a_i)_{\sigma(x)}$. But then $f_{\sigma(x)}((a_i)_{\sigma(x)}) = \ev_x((\id \otimes f)(a)) = 0$ and by injectivity of $f_{\sigma(x)}$, we have $(a_i)_{\sigma(x)} = 0$. Since this holds for each $x \in U_i$, we have $\chi_{U_i} \otimes a_i = 0$ and therefore $a = 0$.

If $n \geq 2$, we may assume by the same disjointification procedure that for each $1 \leq k \leq n$ the $q(U_i) \subseteq O_k$ are pairwise disjoint for $i \in I$. Let $I' \subseteq I$ be the set of $i$ such that $q(U_i) \subseteq O_k$ for some $k < n$. Pick some $j \in I \setminus I'$. Then there is a compact open set $U'_j \subseteq U_j$ such that $\chi_{U_j} \otimes a_j = \chi_{U'_j} \otimes a_j$ and $q(U'_j) \subseteq \bigcup_{i \in I'} q(U_i)$, constructed as follows. For any $u \in U_j$ with $q(u) \notin \bigcup_{i \in I'} q(U_i)$, we may conclude that $a_{\sigma(u)} = 0$ by evaluation at $u$ and injectivity of $f$. Then $a_y = 0$ for $y$ in a compact open neighbourhood $W_u$ of $\sigma(u)$ in $\sigma(U_j)$. Moving $W_u$ for each $u$ to $X/G$ via $U_j$ we obtain a compact open cover of the compact set $q(U_j) \setminus \bigcup_{i \in I'} q(U_i)$ which therefore has a finite subcover. We obtain finitely many compact opens in $U_j$ which we subtract from $U_j$ to obtain $U'_j$. By construction $a$ vanishes on $\sigma(U_j \setminus U_j')$ and $q(U'_j) \subseteq \bigcup_{i \in I'} q(U_i)$.

We now cover $q(U'_j)$ by finitely many compact open subsets $V_1,\dots,V_N \subseteq q(U'_j)$ each contained in $q(U_i)$ for some $i \in I'$, and we may assume that these are disjoint. Thus $\chi_{U'_j} = \sum_{m = 1}^N \chi_{q^{-1}(V_m)\cap U'_j}$, and so we are able to reduce to the case where $a = \sum_{i \in I} \chi_{U_i} \otimes a_i$ such that for each $i \in I$ there is $k < n$ with $q(U_i) \subseteq O_k$.
\end{proof}
\end{prop}

For a flat $G$-module $\Z[Y]$ as in Proposition \ref{flat module result}, the tensor product over $G$ with a left $G$-module of the form $\Z[Z]$ corresponds to the fibre product $Y \times_G Z$ of spaces over $G$, which is the quotient of $Y \times_{G^0} Z$ by the diagonal action of $G$.

\begin{prop}\label{tensor product as fibre product}
Let $G$ be an ample groupoid, let $Y$ be a basic \'etale right $G$-space with anchor map $\sigma \colon Y \to G^0$ and let $Z$ be a totally disconnected left $G$-space. Then $Y \times_G Z$ is totally disconnected and locally LCH, and there is an isomorphism $\kappa \colon \Z[Y] \otimes_G \Z[Z] \cong \Z[Y \times_G Z]$ given on simple tensors by
\begin{align*}
\kappa \colon \Z[Y] \otimes_{G} \Z[Z] & \to \Z[Y \times_G Z] \\
 \xi \otimes \eta & \mapsto  \kappa(\xi, \eta) \\
 & \quad \; \; [y,z]_G \mapsto \sum_{g \in G^{\sigma(y)}} \xi(y \bcdot g) \eta( g^{-1} \bcdot z).
\end{align*}
\begin{proof}
The fibre product $Y \times_{G^0} Z$ is totally disconnected and locally LCH, and the diagonal action $G \acts Y \times_{G^0} Z$ is basic, from which it follows that $Y \times_G Z$ is totally disconnected and locally LCH.

Let $\rho \colon Z \to {G^0}$ be the anchor map for $G \acts Z$, and let $\cal S$ be the set of pairs $(U,V)$ of compact open Hausdorff subsets $U \subseteq Y$ and $V \subseteq Z$ such that $Y \to Y/G$ and $\sigma \colon Y \to {G^0}$ are injective on $U$ and $\rho(V) \subseteq \sigma(U)$. The sets $U \times_{{G^0}} V$ with $(U,V) \in \cal S$ form a basis of compact open Hausdorff sets in $Y \times_{{G^0}} Z$, so their images $q(U \times_{{G^0}} V)$ under the local homeomorphism $q\colon Y \times_{G^0} Z \to Y \times_G Z$ form a basis of compact open Hausdorff sets in $Y \times_G Z$. For $\xi \in \Z[Y]$ and $\eta \in \Z[Z]$ define $\kappa(\xi, \eta) \colon Y \times_G Z \to \Z$ at $[y,z]_G \in Y \times_G Z$ by
\[ \kappa(\xi, \eta) \colon [y,z]_G \mapsto \sum_{g \in G^{\sigma(y)}} \xi(y \bcdot g) \eta( g^{-1} \bcdot z). \]
By construction $\kappa(\chi_U, \chi_V) = \chi_{q(U \times_{G^0} V)}$ for each $(U,V) \in \cal S$ so we obtain a balanced bilinear map $\kappa \colon \Z[Y] \times \Z[Z] \to \Z[Y \times_G Z]$ whose image generates $\Z[Y \times_G Z]$. There is therefore a surjective homomorphism
\begin{align*}
\kappa \colon \Z[Y] \otimes_G \Z[Z] & \to \Z[Y \times_G Z] \\
\xi \otimes \eta & \mapsto \kappa(\xi, \eta)
\end{align*}
Consider the cover $\mathcal R = \{ q(U \times_{G^0} V) \suchthat (U,V) \in \cal S \}$ of compact open Hausdorff sets in $Y \times_G Z$ which is closed under compact open subsets. Using Lemma \ref{basiscompatibility}, we define an inverse $\psi \colon \Z[Y \times_G Z] \to \Z[Y] \otimes_G \Z[Z]$ to $\kappa$ by setting 
\[ \psi(\chi_{q(U \times_{G^0} V)}) := \chi_U \otimes \chi_V. \]
We need to check that this is well-defined and respects disjoint unions within $\cal R$. First, suppose that $(U_1,V_1), (U_2,V_2) \in \cal S$ with $q(U_1 \times_{G^0} V_1) = q(U_2 \times_{G^0} V_2)$. Define
\[ W := \left\{ g \in G \suchthat \text{there are} \; u_1 \in U_1 \; \text{and} \; u_2 \in U_2 \; \text{such that} \; u_1 \bcdot g = u_2 \right\}. \]
The action $Y \rightacts G$ is basic and $\sigma$ restricts to homeomorphisms on $U_1$ and $U_2$, from which it follows that $W$ is a compact open bisection in $G$. Furthermore, because $q(U_1 \times_{G^0} V_1) = q(U_2 \times_{G^0} V_2)$ we obtain that $\rho(V_1) \subseteq r(W)$, $\rho(V_2) \subseteq s(W)$ and $W \bcdot V_2 = V_1$. We may then calculate
\[ \chi_{U_1} \otimes \chi_{V_1} = \chi_{U_1} \otimes \chi_W \bcdot \chi_{V_2} = \chi_{U_1} \bcdot \chi_W \otimes \chi_{V_2} = \chi_{U_2} \bcdot \chi_{r(W)} \otimes \chi_{V_2} = \chi_{U_2} \otimes \chi_{V_2}. \]
Now suppose that $q(U_1 \times_{G^0} V_1) \sqcup q(U_2 \times_{G^0} V_2) = q(U \times_{G^0} V)$. By the above argument we may assume that $U_1 \times_{G^0} V_1 \sqcup U_2 \times_{G^0} V_2 = U \times_{G^0} V$, so $V_1 \sqcup V_2 = V$ and we can write $\chi_U \otimes \chi_V = \chi_U \otimes \chi_{V_1} + \chi_U \otimes \chi_{V_2} = \chi_{U_1} \otimes \chi_{V_1} + \chi_{U_2} \otimes \chi_{V_2}$. By Lemma \ref{basiscompatibility}, $\psi$ extends uniquely to a homomorphism $\psi \colon \Z[Y \times_G Z] \to \Z[Y] \otimes_{G} \Z[Z]$ such that $\kappa \circ \psi = 1$ by construction. The elements $\chi_U \otimes \chi_V$ for $(U,V) \in \cal S$ generate $\Z[Y] \otimes_{G} \Z[Z]$, so therefore $\psi$ is an inverse to $\kappa$.
\end{proof}
\end{prop}
Proposition \ref{tensor product as fibre product} extends \cite[Lemma 2.5]{BDGW23} which describes the coinvariants $\Z[X]_G$ of the module of a free, proper, \'etale $G$-space $X$.

\begin{defn}[Coinvariants functor]
Let $G$ be an ample groupoid and let $M$ be a $G$-module. The \textit{coinvariants} $M_G$ of $M$ is the abelian group
\[ M_G := \Z[G^0] \otimes_{G} M. \]
This gives us a functor $\Coinv_G \colon G\cat{-Mod} \to \Ab$. We often think of $M_G$ as a quotient of $M$ via the surjective homomorphism $\pi_G \colon m \mapsto [m] \colon M \to M_G$ which sends $m$ to $e \otimes m$ for any idempotent $e \in \Z[G^0]$ such that $e \bcdot m = m$. The kernel of $\pi_G$ is generated by elements of the form $\chi_U \bcdot m - \chi_{s(U)} \bcdot m$ for compact open bisections $U \subseteq G$ and elements $m \in M$.
\end{defn}

To discuss derived functors, we recall the fundamental lemma of homological algebra \cite[Comparison Theorem 6.16]{Rotman09}.
\begin{lem}[Fundamental lemma of homological algebra]\label{fundamental lemma of homological algebra}
Let $R$ be a locally unital ring, let $f \colon A \to B$ be a map of $R$-modules, let $P_\bullet \to A$ be a projective resolution and $Q_\bullet \to B$ a resolution. Then there is a chain map $P_\bullet \to Q_\bullet$ over $f$ which is unique up to chain homotopy.
\end{lem}

\begin{defn}[Derived functors]
Given an additive functor $F \colon \lMod R \to \Ab$, the derived functors $\mathbb L_n F \colon \lMod R \to \Ab$ for $n \in \bb N$ are defined as follows. For each $R$-module $A$, there is a projective resolution $P_\bullet \to A$ because $\lMod R$ has enough projectives. We set $\bb L_n F (A) = \H_n(F(P_\bullet))$, which is independent of the choice of projective resolutions up to canonical isomorphisms by the fundamental lemma of homological algebra. Given a morphism $f \colon A \to B$ of $R$-modules and a projective resolution $Q_\bullet \to B$, there is a chain map $P_\bullet \to Q_\bullet$ over $f$ which is unique up to chain homotopy, and therefore for each $n \in \bb N$ a unique induced map $\H_n(F(P_\bullet)) \to \H_n(F(Q_\bullet))$, which we set $\bb L_n (f) \colon \bb L_n(A) \to \bb L_n(B)$ to be.
\end{defn}
The Tor group $\Tor^R_n(A,B)$ of a right $R$-module $A$ and a left $R$-module $B$ is the $n$th left derived functor of $A \otimes_R -$ at $B$. We use this to define groupoid homology\footnote{This definition is shown to be equivalent to the usual definition under mild assumptions in \cite{BDGW23} and in full generality in \cite{Li22}.}.
\begin{defn}[Groupoid homology]
Let $G$ be an ample groupoid and let $M$ be a $G$-module. The \textit{groupoid homology $\H_n(G;M)$ of $G$ with coefficients in $M$} is the $n$th $\Tor$ group $\Tor^G_n(\Z [G^0], M)$. This may be identified with $\bb L_n \Coinv_G (M)$.
\end{defn}

The Tor groups of a locally unital ring may be computed with flat resolutions in place of projective resolutions, as in the case of unital rings \cite[Corollary 10.23]{Rotman09}. This may be justified through use of the multiplier ring, see Remark \ref{multiplier ring}.

\begin{example}[Bar resolution]\label{bar resolution example}
Let $G$ be an ample groupoid. There is an explicit flat resolution $(\Z[G^{\bullet+1}], \partial_\bullet)$ of the (left) $G$-module $\Z[G^0]$ 
\begin{equation}\label{bar resolution}
  \cdots \xrightarrow{\partial_{n+1}} \Z[G^{n+1}] \xrightarrow{\partial_n} \cdots \xrightarrow{\partial_2} \Z[G^2] \xrightarrow{\partial_1} \Z[G^1] \xrightarrow{\partial_0} \Z[G^0] \to 0 
\end{equation}
called the \textit{bar resolution}. For $n \geq 0$ we consider the space $G^{n+1}$ of composable $n+1$-tuples as a left $G$-space, whose $G$-module $\Z[G^{n+1}]$ is flat. We set $\partial_0 := s_* \colon \Z[G^1] \to \Z[G^0]$. For $n \geq 1$ and $0 \leq i \leq n$ we define face maps $\partial^n_i \colon G^{n+1} \to G^{n}$ by
\begin{equation*} \partial^n_i \colon (g_0,\dots,g_n) \mapsto 
\begin{cases}[c] 
(g_0,\dots,g_i g_{i+1},\dots,g_n) & \text{if} \; i < n, \\ 
(g_0,\dots,g_{n-1}) & \text{if} \; i = n. 
\end{cases}  
\end{equation*}
The face maps are $G$-equivariant local homeomorphisms and therefore induce $G$-module maps $(\partial^n_i)_* \colon \Z[G^{n+1}] \to \Z[G^n]$. The boundary maps $\partial_n \colon \Z[G^{n+1}] \to \Z[G^n]$ are given for $n \geq 1$ by
\[ \partial_n := \sum_{i=0}^n (-1)^i (\partial^n_i)_*. \]
The exactness of the bar resolution (\ref{bar resolution}) is witnessed by a chain homotopy induced by local homeomorphisms $h_n \colon G^n \to G^{n+1}$. These are defined for $n \geq 1$ by
\begin{equation*}
h_n \colon (g_0,\dots,g_{n-1})  \mapsto (r(g_0),g_0,\dots,g_{n-1})  
\end{equation*}
and $h_0$ is the inclusion $G^0 \subseteq G^1$. Taking the coinvariants of the bar resolution, we obtain the chain complex
\begin{equation}\label{concrete chain complex for homology}
  \cdots \xrightarrow{(\partial_{n+1})_G} \Z[G^{n}] \xrightarrow{(\partial_n)_G} \cdots \xrightarrow{(\partial_2)_G} \Z[G^1] \xrightarrow{(\partial_1)_G} \Z[G^0]  \to 0.
\end{equation}
The boundary maps are given by $(\partial_n)_G = \sum_{i=0}^n (-1)^i (\epsilon^n_i)_*$, where for $n \geq 1$, the face map $\epsilon^n_i \colon G^n \to G^{n-1}$ is defined by
\begin{equation*}
\epsilon^n_i \colon (g_1,\dots,g_n) \mapsto 
\begin{cases}[c]
(g_2,  \dots, g_n) & \text{if} \; i = 0, \\
(g_1, \dots, g_i g_{i+1}, \dots, g_n) & \text{if} \; 0 < i < n, \\
(g_1, \dots, g_{n-1}) & \text{if} \; i = n.
\end{cases}
\end{equation*}
The chain complex (\ref{concrete chain complex for homology}) is typically used as a concrete way to define $\H_*(G)$ (see \cite[Section 3.1]{Matui12}). 
\end{example}

The homology of an ample groupoid $G$ can be computed with even more general types of resolutions called left $G$-acyclic resolutions.

\begin{defn}
Let $G$ be an ample groupoid. A $G$-module $A$ is \textit{left $G$-acyclic} if $\H_n(G;A) = 0$ for each $n \geq 1$.
\end{defn}

\begin{prop}
Let $G$ be an ample groupoid, let $M$ be a $G$-module, let $P_\bullet \to M$ be a projective resolution and let $Q_\bullet \to M$ be a left $G$-acyclic resolution. Then the unique (up to homotopy) chain map $P_\bullet \to Q_\bullet$ induces an isomorphism $\H_n(G;M) \cong \H_n((Q_\bullet)_G)$ for each $n$. Furthermore, given a $G$-module map $f \colon M \to N$, left $G$-acyclic resolutions $Q_\bullet \to M$ and $Q'_\bullet \to N$ and a chain map $Q_\bullet \to Q'_\bullet$ over $f$, the induced map in homology $\H_n((Q_\bullet)_G) \to \H_n((Q'_\bullet)_G)$ may be identified with $\H_n(G;f) \colon \H_n(G;M) \to \H_n(G;N)$.
\begin{proof}
The first part of the statement is covered in the setting of unital rings in \cite[2.4.3]{Weibel94}, including the linked exercise, but holds in general for abelian categories with enough projectives (see \cite[Section 13.3]{KasSch06}). Alternatively, viewing $G$-modules as modules over the multiplier ring $M(\Z[G])$, we preserve exactness, projectivity and tensor products, and therefore also $\Tor$ groups. The isomorphism $\H_n(G;M) \cong \H_n((Q_\bullet)_G)$ can be deduced by working over the unital ring $M(\Z[G])$. The identification of the induced map in homology then follows from the fundamental lemma of homological algebra.
\end{proof}
\end{prop}

\begin{defn}
Let $G$ be an ample groupoid and let $H$ be an open subgroupoid of $G$. The inclusion $\Z[H] \subseteq \Z[G]$ gives rise to the \textit{subgroupoid restriction functor} $\Res^H_G \colon \lMod G \to \lMod H$ which for a $G$-module $M$ returns the $H$-module
\[ \Res^H_G  M = \Z[H] \bcdot M.\]
The \textit{subgroupoid induction functor} $\Ind^G_H \colon \lMod H \to \lMod G$ is given by the tensor product with the bimodule $\Z[G_{H^0}]$, so for an $H$-module $N$ we have
\[ \Ind^G_H N = \Z[G_{H^0}] \otimes_H N. \]
\end{defn}
The restriction functor $\Res^H_G$ sends a $G$-module map to its restriction and is therefore exact. The action $G_{H^0} \rightacts H$ is basic and \'etale, so $\Ind^G_H$ is an exact functor by Proposition \ref{flat module result}.

\begin{prop}[Induction-restriction adjunction]\label{induction restriction}
Let $G$ be an ample groupoid and let $H \subseteq G$ be an open subgroupoid. Then there is an adjunction
\[ (\eta, \epsilon) \colon \Ind^G_H \dashv \Res^H_G \]
with unit and counit
\begin{align*}
\eta \colon \id_{\lMod H} & \Rightarrow \Res^H_G \Ind^G_H \\
\epsilon \colon \Ind^G_H \Res^H_G & \Rightarrow \id_{\lMod G}  
\end{align*}
such that for $M \in \lMod H$, $N \in \lMod G$, $\zeta \in \Z[H]$, $m \in M$, $\xi \in \Z[G_{H^0}] \subseteq \Z[G]$ and $n \in \Z[H] \bcdot N$, we have
\begin{align*}
\eta_M \colon M & \to \Z[H] \bcdot (\Z[G_{H^0}] \otimes_H M) \\
\zeta \bcdot m & \mapsto    \zeta \otimes m,   \\
\epsilon_N \colon \Z[G_{H^0}] \otimes_H ( \Z[H] \bcdot N) & \to N \\
\xi \otimes n & \mapsto \xi \bcdot n.
\end{align*}
\begin{proof}
For each $H$-module $M$, well-definition of the map $\eta_M$ follows from the application of a local unit on the left of each $\zeta \in \Z[H]$. This is an $H$-module map and natural in $M$ by construction. For each $G$-module $N$ the above formula describes a balanced bilinear map through which we obtain $\epsilon_N$, which is a $G$-module map and natural in $N$. To verify the counit-unit equations we must check that
\begin{align*}
  \epsilon_{\Ind^G_H M} \circ \Ind^G_H(\eta_M)   & = \id_{\Ind^G_H M}, \\
  \Res^H_G (\epsilon_N) \circ \eta_{\Res^H_G N}  & = \id_{\Res^H_G N}.
\end{align*}
Let $m \in M$, $n \in N$, $\xi \in \Z[G_{H^0}]$ and $\zeta \in \Z[H]$. We then calculate
\begin{align*}
\epsilon_{\Ind^G_H M} (\Ind^G_H (\eta_M)(\xi \otimes \zeta \bcdot m)) & = \epsilon_{\Ind^G_H M} (\xi \otimes \eta_M(\zeta \bcdot m)) \\
& = \epsilon_{\Ind^G_H M} (\xi \otimes \zeta \otimes m) \\
& = \xi \bcdot \zeta \otimes m \\
& = \xi \otimes \zeta \bcdot m.
\end{align*}
Letting $e \in \Z[H]$ be a left unit for $\zeta$, we verify the second equation:
\begin{align*}
\Res^H_G(\epsilon_N)(    \eta_{\Res^H_G N} ( \zeta \bcdot n )  ) & = \Res^H_G(\epsilon_N)(    \eta_{\Res^H_G N} ( e \bcdot  \zeta \bcdot n )  ) \\
& = \Res^H_G(\epsilon_N)(       e \otimes \zeta \bcdot n     ) \\
& = \epsilon_N ( e \otimes \zeta \bcdot n ) \\
& = \zeta \bcdot n
\end{align*}
\end{proof}
\end{prop}

\begin{lem}[Shapiro's Lemma]\label{Shapiros lemma}
Let $G$ be an ample groupoid and let $H \subseteq G$ be an open subgroupoid of $G$. Then for any $H$-module $M$, there is an isomorphism 
\[\H_*(G;\Ind^G_H M) \cong \H_*(H;M).\]
\begin{proof}
Let $F_\bullet \to M$ be a flat resolution in $\lMod H$. Then $\Ind^G_H F_\bullet \to \Ind^G_H M$ is exact by exactness of $\Ind^G_H$. The $G$-module $\Z[G_{H^0}]$ is flat and it follows that $\Ind^G_H F_\bullet \to \Ind^G_H M$ is a flat resolution in $\lMod G$. By Proposition \ref{tensor product as fibre product} there is an isomorphism $\Z[G^0] \otimes_G \Z[G_{H^0}] \cong \Z[H^0]$ of right $H$-modules. Through this there is a natural isomorphism $(\Ind^G_H N)_G \cong N_H$ for each $H$-module $N$, and so a chain isomorphism $(\Ind^G_H F_\bullet)_G \cong (P_\bullet)_H$. This witnesses the claimed isomorphism in homology.
\end{proof}
\end{lem}

\begin{cor}\label{Shapiro corollary}
Let $G$ be an ample groupoid with unit space $X$ and let $M$ be an $X$-module. Then $\Ind^G_{X} M$ is left $G$-acyclic.
\end{cor}

\begin{example}[Bar resolution of a module]
Let $G$ be an ample groupoid with unit space $X$, let $M$ be a $G$-module and consider the functor 
\[L = \Ind^G_X \Res^X_G \colon \lMod G \to \lMod G.\]
From the counit $\epsilon \colon L \Rightarrow \id$ of the adjunction $\Ind^G_X \dashv \Res^X_G$ we may construct a chain complex
\begin{equation}\label{bar resolution of a module}
 \cdots \to L^{n+1} M \xrightarrow{\delta_n} \cdots \xrightarrow{\delta_1} L M \xrightarrow{\pi_0} M \to 0 
\end{equation}
with $\pi_0 = \epsilon_M \colon LM \to M$ and for $n \geq 1$
\[ \delta_n = \sum_{i=0}^n (-L)^i \epsilon_{L^{n-i} M} \colon L^{n+1} M \to L^n M. \]
This is exact because applying $\Res^X_G$ does not change the underlying chain complex of abelian groups, and in $\lMod X$ there is a contracting homotopy coming from the unit $\eta \colon \id \Rightarrow \Res^X_G \Ind^G_X$ of the adjunction (see \cite[8.6.10, 8.6.11]{Weibel94}). This is given at $n \geq 0$ by $h_n = \eta_{\Res^X_G L^n M} \colon \Res^X_G L^n M \to \Res^X_G L^{n+1} M$. The resolution $L^{\bullet +1}M \to M$ in (\ref{bar resolution of a module}) is the \textit{bar resolution of $M$}. It is left $G$-acyclic by Corollary \ref{Shapiro corollary} and may therefore be used to compute $\H_*(G;M)$. We note that the chain complex $\Coinv_G (L^{\bullet+1} M) = \Z [G^{\bullet}] \otimes_X M$ coincides with the chain complex defining $\H_*(G;M)$ in \cite{ProYam22}. For $M = \Z[X]$ we recover the bar resolution (Example \ref{bar resolution example}).
\end{example}

\section{\'Etale correspondences of ample groupoids}

We refer the reader to \cite{AKM22} and \cite{Miller23a} for an introduction to and examples of \'etale correspondences of \'etale groupoids.

\begin{defn}[Definition 3.1 in \cite{AKM22}]
Let $G$ and $H$ be \'etale groupoids. An \textit{\'etale correspondence} $\Omega \colon G \to H$ is a space $\Omega$ with a left $G$-action and a right $H$-action with anchor maps $\rho\colon \Omega \to G^0$ and $\sigma \colon  \Omega \to H^0$ called the \textit{range} and \textit{source} such that:
\begin{itemize}
\item The $G$-action commutes with the $H$-action, i.e. $\Omega$ is a \textit{$G$-$H$-bispace}.
\item The right action $\Omega \rightacts H$ is free, proper and \'etale.
\end{itemize}
We say that $\Omega \colon G \to H$ is \textit{proper} if the map $\overline{\rho} \colon \Omega/H \to G^0$ induced by $\rho$ is proper. If we want to highlight the correspondence $\Omega$, we may write $\sigma_\Omega$ and $\rho_\Omega$ instead of $\sigma$ and $\rho$. For $x \in G^0$ and $y \in H^0$, we write $\Omega^x$ and $\Omega_y$ for the range and source fibres $\rho^{-1}(x)$ and $\sigma^{-1}(y)$. An \textit{isomorphism} of \'etale correspondences $\Omega, \Lambda \colon G \to H$ is a $G$-$H$-equivariant homeomorphism $\Omega \cong \Lambda$.

The \textit{composition} $\Lambda \circ \Omega$ of \'etale correspondences $\Omega \colon G \to H$ and $\Lambda \colon H \to K$ is given by the space $\Omega \times_H \Lambda$, with $G$ acting through its action on $\Omega$ and $K$ acting through its action on $\Lambda$. This forms an associative composition up to isomorphism, with identities up to isomorphism given at $G$ by the $G$-$G$-bispace $G$.
\end{defn}

If $\Omega \colon G \to H$ is an \'etale correspondence of ample groupoids, then $\Omega$ is a totally disconnected locally LCH space. The space $\Z[\Omega]$ becomes a $G$-$H$-bimodule.

\begin{defn}
Let $\Omega \colon G \to H$ be an \'etale correspondence of ample groupoids. The \textit{induction functor} $\Ind_\Omega \colon \lMod H \to \lMod G$ is given by the tensor product $\Z[\Omega] \otimes_H -$.
\end{defn}

Proposition \ref{flat module result} shows that $\Ind_\Omega$ is an exact functor and Proposition \ref{tensor product as fibre product} gives an explicit description of $\Ind_\Omega$ on $H$-modules of the form $\Z[Y]$ for an $H$-space $Y$ as $\Z[\Omega \times_H Y]$. In particular we have $\Ind_\Omega \Z[H^0] \cong \Z[\Omega/H]$. For \'etale correspondences $\Omega \colon G \to H$ and $\Lambda \colon H \to K$ of ample groupoids, we obtain an isomorphism $\Z[\Omega \times_H \Lambda] \cong \Z[\Omega] \otimes_H \Z[\Lambda]$ and therefore a natural isomorphism $\Ind_{\Lambda \circ \Omega} \cong \Ind_\Omega \Ind_\Lambda$. 

\begin{prop}\label{coinvariantmodulemap}
Let $\Omega \colon G \to H$ be an \'etale correspondence of ample groupoids. There is a map $\delta_\Omega \colon \Z[\Omega]_G  \to \Z[H^0]$ of right $H$-modules such that the following diagram commutes.
\begin{equation*} 
\begin{tikzcd}
\Z[\Omega] \arrow[r, "\sigma_*"] \arrow[d, "\pi_G"'] & \Z[H^0] \\
\Z[\Omega]_G  \arrow[ur, "\delta_\Omega"] &
\end{tikzcd}
\end{equation*}
The map $\delta_\Omega$ satisfies $\delta_\Omega(\eta \otimes \xi) = \sigma_*(\eta \bcdot \xi)$ for each $\eta \in \Z[G^0]$ and $\xi \in \Z[\Omega]$.
\begin{proof}
The balancedness of the bilinear map 
\begin{align*}
 \Z[G^0] \times \Z[\Omega] & \to \Z[H^0]  \\
 (\eta, \xi) & \mapsto \sigma_*(\eta \bcdot \xi)
\end{align*}
follows from the $G$-invariance of $\sigma \colon \Omega \to H^0$. We therefore obtain our desired map $\delta_\Omega \colon \Z[\Omega]_G \to \Z[H^0]$. The diagram commutes because for each $\eta \in \Z[G^0]$ and $\xi \in \Z[\Omega]$, we have $\pi_G(\eta \bcdot \xi) = \eta \otimes \xi$.
\end{proof}
\end{prop}
Through the tensor product over $H$, the map $\delta_\Omega \colon \Z[\Omega]_G \to \Z[H^0]$ induces a map $\delta_\Omega \otimes \id \colon (\Ind_\Omega M)_G \to M_H$ for each $H$-module $M$, and this is natural in $M$. We have compatibility with composition of correspondences:
\begin{prop}\label{coinvariant module map and composition}
Let $\Omega \colon G \to H$ and $\Lambda \colon H \to K$ be \'etale correspondences of ample groupoids and let $\kappa \colon \Z[\Omega] \otimes_H \Z[\Lambda] \cong \Z[\Omega \times_H \Lambda]$ be the isomorphism in Proposition \ref{tensor product as fibre product}. Then the following diagram commutes.
\[ \begin{tikzcd}
{\Z[\Omega]_G \otimes_H \Z[\Lambda]} \arrow[r, "\delta_\Omega \otimes \id"] \arrow[d, "\kappa_G"'] & {\Z[\Lambda]_H} \arrow[d, "\delta_\Lambda"] \\
{\Z[\Omega \times_H \Lambda]_G} \arrow[r, "\delta_{\Lambda \circ \Omega}"]                         & {\Z[K^0]}                                  
\end{tikzcd} \]
\begin{proof}
Let $\eta \in \Z[G^0]$, $\xi \in \Z[\Omega]$ and $\nu \in \Z[\Lambda]$. The simple tensor $\eta \otimes \xi \otimes \nu \in \Z[G^0] \otimes_G \Z[\Omega] \otimes_H \Z[\Lambda]$ is sent to the following element of $\Z[K^0]$ through both routes round the diagram.
\begin{align*}
K^0 & \to \Z  \\
 z & \mapsto \sum_{\lambda \in \Lambda_z} \sum_{\omega \in \Omega_{\rho(\lambda)}} \eta(\rho(\omega))\xi(\omega)\nu(\lambda)  
\end{align*}
\end{proof}
\end{prop}

\begin{thm}\label{inducedmapinhomology}
Let $\Omega \colon G \to H$ be an \'etale correspondence of ample groupoids, let $M$ be a $G$-module, $N$ an $H$-module and $f \colon M \to \Ind_\Omega N$ a $G$-module map. Then there is an induced map in homology
\[ \H_*(\Omega, f) \colon \H_*(G;M) \to \H_*(H; N)\]
such that for any
\begin{itemize}
\item left $G$-acyclic resolution $P_\bullet \to M$,
\item left $H$-acyclic resolution $Q_\bullet \to N$,
\item and chain map $\tilde f \colon P_\bullet \to \Ind_\Omega Q_\bullet$ lifting $f \colon M \to \Ind_\Omega N$,
\end{itemize}
the chain map $(\delta_\Omega \otimes \id) \circ \tilde f_G \colon (P_\bullet)_G \to (Q_\bullet)_H$ shown below induces $\H_{*}(\Omega, f)$ in homology.
\[\begin{tikzcd}
	\cdots & {(P_n)_G} & \cdots & {(P_0)_G} \\
	\cdots & {(\Ind_\Omega Q_n)_G} & \cdots & {(\Ind_\Omega Q_0)_G} \\
	\cdots & {(Q_n)_H} & \cdots & {(Q_0)_H} 
	\arrow[from=1-2, to=2-2, "(\tilde f_n)_G"]
	\arrow[from=1-2, to=1-3]
	\arrow[from=1-3, to=1-4]
	\arrow[from=2-2, to=2-3]
	\arrow[from=2-3, to=2-4]
	\arrow[from=1-4, to=2-4, "(\tilde f_0)_G"]
	\arrow[from=2-1, to=2-2]
	\arrow[from=1-1, to=1-2]
	\arrow["{\delta_\Omega \otimes \id}", from=2-2, to=3-2]
	\arrow[from=3-2, to=3-3]
	\arrow[from=3-1, to=3-2]
	\arrow[from=3-3, to=3-4]
	\arrow["{\delta_\Omega \otimes \id}", from=2-4, to=3-4]
\end{tikzcd}\]
\begin{proof}
To construct the map $\H_{*}(\Omega, f)$, we may consider arbitrary projective resolutions $P'_\bullet \to M$ and $Q'_\bullet \to N$. We obtain a resolution $\Ind_\Omega Q'_\bullet \to \Ind_\Omega N$. A chain map $\tilde f' \colon P'_\bullet \to \Ind_\Omega Q'_\bullet$ lifting $f \colon M \to \Ind_\Omega N$ exists by the fundamental lemma of homological algebra (Lemma \ref{fundamental lemma of homological algebra}). We may then set $\H_{*}(\Omega,f)$ to be $\H_*((\delta_\Omega \otimes \id) \circ \tilde f'_G) \colon \H_*((P'_\bullet)_G) \to \H_*((Q'_\bullet)_H)$ after identifying these with $\H_*(G;M)$ and $\H_*(H;N)$.

Now given $P_\bullet \to M$, $Q_\bullet \to N$ and $\tilde f \colon P_\bullet \to \Ind_\Omega Q_\bullet$ as in the statement of the theorem, there are chain maps $\pi \colon P_\bullet' \to P_\bullet$ and $\tau \colon Q'_\bullet \to Q_\bullet$. By the fundamental lemma of homological algebra, $\tilde f \circ \pi$ is chain homotopic to $\Ind_\Omega \tau \circ \tilde f'$. The following diagram of chain complexes therefore commutes up to chain homotopy. 
\[\begin{tikzcd}[column sep = 4em]
(P'_\bullet)_G \arrow[r, "\pi_G"] \arrow[d, "\tilde f'_G"] & (P_\bullet)_G \arrow[d, "\tilde f_G"] \\
(\Ind_\Omega Q'_\bullet)_G \arrow[r, "(\Ind_\Omega \tau )_G"] \arrow[d, "\delta_\Omega \otimes \id"] & (\Ind_\Omega Q_\bullet)_G \arrow[d,"\delta_\Omega \otimes \id"] \\
(Q'_\bullet)_H \arrow[r, "\tau_H"] & (Q_\bullet)_H
\end{tikzcd} \]
We deduce that $(\delta_\Omega \otimes \id) \circ \tilde f'_G$ and $(\delta_\Omega \otimes \id) \circ \tilde f_G$ induce the same map in homology.
\end{proof}
\end{thm}
In particular, if $\Omega \colon G \to H$ is a proper correspondence, then there is a $G$-module map $\orho^* \colon \Z[G^0] \to \Ind_\Omega \Z[H^0]$ induced by the proper $G$-equivariant map $\orho \colon \Omega/H \to G^0$.

\begin{cor}\label{main corollary}
Let $G$ and $H$ be ample groupoids. Then any proper \'etale correspondence $\Omega \colon G \to H$ induces a map
\[ \H_*(\Omega) \colon \H_*(G) \to \H_*(H) \]
in homology. This has the property that for any left $G$-acyclic resolution $P_\bullet \to \Z[G^0]$, left $H$-acyclic resolution $Q_\bullet \to \Z[H^0]$ and chain map $f \colon P_\bullet \to \Ind_\Omega Q_\bullet$ over $\orho^* \colon \Z[G^0] \to \Ind_\Omega \Z[H^0]$, $\H_*(\Omega)$ is induced by the chain map
\[ (\delta_\Omega \otimes \id) \circ f_G \colon (P_\bullet)_G \to (Q_\bullet)_H. \]
\end{cor}

\begin{prop}
Let $\Omega \colon G \to H$ and $\Lambda \colon H \to K$ be proper \'etale correspondences of ample groupoids. Then we have 
\[\H_*(\Lambda \circ \Omega) = \H_*(\Lambda) \circ \H_*(\Omega) \colon \H_*(G) \to \H_*(K). \]
\begin{proof}
Consider the proper maps $\orho \colon (\Omega \times_H \Lambda)/K \to G^0$, $\orho_\Omega \colon \Omega/H \to G^0$ and $\orho_\Lambda \colon \Lambda/K \to H^0$. Let $P_\bullet \to \Z[G^0]$, $Q_\bullet \to \Z[H^0]$ and $R_\bullet \to \Z[K^0]$ be left $G$, $H$ and $K$-acyclic resolutions, and let $f \colon P_\bullet \to \Ind_\Omega Q_\bullet$ and $g \colon Q_\bullet \to \Ind_\Lambda R_\bullet$ be chain maps over $\orho_\Omega^* \colon \Z[G^0] \to \Ind_\Omega \Z[H^0]$ and $\orho_\Lambda^* \colon \Z[H^0] \to \Ind_\Lambda \Z[K^0]$. Then
\[ \begin{tikzcd}
P_\bullet \arrow[r, "f"] & \Ind_\Omega Q_\bullet \arrow[r, "\id \otimes g"] & \Ind_\Omega \Ind_\Lambda R_\bullet \arrow[r, "\kappa \otimes \id"] & \Ind_{\Lambda \circ \Omega} R_\bullet
\end{tikzcd} \]
is a chain map over
\[ \begin{tikzcd}
{\Z[G^0]} \arrow[r, "\orho_\Omega^*"] & {\Ind_\Omega \Z[H^0]} \arrow[r, "\id \otimes \orho_\Lambda^*"] & {\Ind_\Omega \Ind_\Lambda \Z[K^0]} \arrow[r, "\kappa \otimes \id"] & {\Ind_{\Lambda \circ \Omega} \Z[K^0],}
\end{tikzcd} \]
where $\kappa \colon \Z[\Omega] \otimes_H \Z[\Lambda] \cong \Z[\Omega \times_H \Lambda]$ is the isomorphism in Proposition \ref{tensor product as fibre product}. Under the identifications of modules induced from $H$-space modules and $K$-space modules with quotient $G$-space modules, these maps are induced by the proper maps
\[ \begin{tikzcd}
G^0 & \Omega/H \arrow[l, "\orho_\Omega"] & \Omega \times_H (\Lambda/K) \arrow[l, "\id \times \orho_\Lambda"] & (\Omega \times_H \Lambda)/K. \arrow[l, "\cong"]
\end{tikzcd} \] 
The composition of these maps is $\orho \colon (\Omega \times_H \Lambda)/K \to G^0$, and we may conclude by Corollary \ref{main corollary} that $\H_*(\Lambda \circ \Omega)$ is induced by the chain map 
\[  (\delta_{\Lambda \circ \Omega} \otimes \id) \circ ((\kappa \otimes \id) \circ (\id \otimes g) \circ f)_G \colon (P_\bullet)_G \to (R_\bullet)_K. \]
Consider the following diagram.
\[ \begin{tikzcd}[column sep = large]
(P_\bullet)_G \arrow[r, "f_G"] & (\Ind_\Omega Q_\bullet)_G \arrow[r, "(\id \otimes g)_G"] \arrow[d, "\delta_\Omega \otimes \id"] & (\Ind_\Omega \Ind_\Lambda R_\bullet)_G \arrow[r, "(\kappa \otimes \id)_G"] \arrow[d, "\delta_\Omega \otimes \id"] & (\Ind_{\Lambda \circ \Omega} R_\bullet )_G \arrow[d, "\delta_{\Lambda \circ \Omega} \otimes \id"] \\
                               & (Q_\bullet)_H \arrow[r, "g_H"]                                                                  & (\Ind_\Lambda R_\bullet)_H \arrow[r, "\delta_\Lambda \otimes \id"]                                                & (R_\bullet)_K                                                                                   
\end{tikzcd} \]
By Corollary \ref{main corollary}, the lower route induces $\H_*(\Lambda) \circ \H_*(\Omega)$, so it suffices to verify that the diagram commutes. The right square commutes by Proposition \ref{coinvariant module map and composition}, and the left square commutes by naturality of $\delta_\Omega \otimes \id \colon (\Ind_\Omega M)_G \to M_H$ in $H$-modules $M$.
\end{proof}
\end{prop}

\begin{example}[Induced map in homology from an \'etale homomorphism]\label{map in homology from etale homomorphism}
Let $\varphi \colon G \to H$ be an \'etale homomorphism of ample groupoids. Recall from Example \ref{bar resolution} that the chain complex $(\Z[G^\bullet], (\partial_\bullet)_G)$ (\ref{concrete chain complex for homology}) computes the homology of $G$. For each $n \geq 0$ the local homeomorphism $\varphi^{n} \colon G^{n} \to H^{n}$ induces a homomorphism $\varphi^{n}_* \colon \Z[G^{n}] \to \Z[H^{n}]$. These form a chain map which induces a map in homology $\H_*(\varphi) \colon \H_*(G) \to \H_*(H)$.

The associated correspondence $\Omega_\varphi \colon G \to H$ is the space $G^0 \times_{H^0} H$. Consider the bar resolutions $\Z[G^{\bullet+1}] \to \Z[G^0]$ and $\Z[H^{\bullet+1}] \to \Z[H^0]$. For each $n \geq 0$ the induced module $\Ind_{\Omega_\varphi} \Z[H^{n+1}]$ has underlying abelian group $\Z[G^0 \times_{H^0} H^{n+1}]$, with the $G$-module structure coming from the action $G \acts G^0 \times_{H^0} H^{n+1}$. Consider the following local homeomorphisms.
\begin{align*}
\psi_n \colon G^{n+1} & \to G^0 \times_{H^0} H^{n+1} & \mu_n \colon G^0 \times_{H^0} H^{n+1} & \to H^n \\
(g_0, \dots, g_n) & \mapsto (r(g),\varphi(g_0),\dots,\varphi(g_n)) & (x,h_0,\dots,h_n) & \mapsto (h_1,\dots,h_n) \\
\epsilon_{G,n} \colon G^{n+1} & \to G^n & \epsilon_{H,n} \colon H^{n+1} & \to H^n \\
(g_0,\dots,g_n) & \mapsto (g_1,\dots,g_n) & (h_0,\dots,h_n) & \mapsto (h_1,\dots,h_n)
\end{align*}
The map $\psi_n$ is $G$-equivariant, while $\mu_n$ and $\epsilon_{G,n}$ are $G$-invariant and $\epsilon_{H,n}$ is $H$-invariant. We obtain a chain map 
\begin{equation*}
f := (\psi_\bullet)_* \colon \Z[G^{\bullet +1}] \to \Ind_{\Omega_\varphi} \Z[H^{\bullet +1}]
\end{equation*}
over the identity $\Z[G^0] \to \Z[G^0] = \Ind_{\Omega_{\varphi}} \Z[H^0]$. By Corollary \ref{main corollary}, $H_*(\Omega_\varphi)$ is induced by the chain map $( \delta_{\Omega_\varphi} \otimes \id ) \circ f_G \colon \Z[G^{\bullet+1}]_G \to \Z[H^{\bullet+1}]_H$. The coinvariants $\Z[G^{n+1}]_G$ and $\Z[H^{n+1}]_H$ are given by $\Z[G^n]$ and $\Z[H^n]$ with quotient maps induced by $\epsilon_{G,n}$ and $\epsilon_{H,n}$. Under these identifications, $\delta_{\Omega_\varphi} \otimes \id  \colon (\Ind_{\Omega_{\varphi}} \Z[H^{n+1}])_G \to \Z[H^{n+1}]_H$ is induced by the $G$-invariant map $\mu_n$. 
\begin{equation*}
\begin{tikzcd}[column sep = large]
{\Z[G^{n+1}]} \arrow[r, "(\psi_n)_*"] \arrow[d, "(\epsilon_{G,n})_*", "\pi_G"']                       & {\Z[G^0 \times_{H^0} H^{n+1}]} \arrow[d, "\pi_G"'] \arrow[rdd, "(\mu_n)_*", bend left] &           \\
{\Z[G^n]} \arrow[r, "f_G"] \arrow[rrd, "(\varphi^n)_*"'] & {(\Z[G^0 \times_{H^0} H^{n+1}])_G} \arrow[rd, "\delta_{\Omega_\varphi} \otimes \id"]  &           \\
                                                                                &                                                                                        & {\Z[H^n]}
\end{tikzcd}
\end{equation*}
The equality $\mu_n \circ \psi_n = \varphi^n \circ \epsilon_{G,n} \colon G^{n+1} \to H^n$ of local homeomorphisms implies that $(\delta_{\Omega_\varphi} \otimes \id) \circ ((\psi_n)_*)_G = \varphi^n_* \colon \Z[G^n] \to \Z[H^n]$. It follows that $\H_*(\Omega_{\varphi}) = \H_*(\varphi)$, so Corollary \ref{main corollary} recovers the standard functoriality of groupoid homology with respect to \'etale homomorphisms.
\end{example}

\begin{example}[Induced map in homology from an action correspondence]\label{map in homology from action correspondence}
Let $G$ be an ample groupoid, let $X$ be a totally disconnected LCH $G$-space with a proper anchor map $\tau \colon X \to G^0$ and let $H = G \ltimes X$ be the action groupoid. For each $n$ consider the map $\tau_{n} \colon H^n \to G^n$ forgetting the elements of $X$. For each $U \subseteq G^n$ in an open Hausdorff cover of $G^n$ coming from bisections in $G$, the preimage $\tau_n^{-1}(U)$ is Hausdorff and the restriction $\tau_n \colon \tau_n^{-1}(U) \to U$ is proper, and we therefore obtain a map $\tau_{n}^* \colon \Z[G^n] \to \Z[H^n]$. These form a chain map which induces a map in homology $\H_*(\tau) \colon \H_*(G) \to \H_*(H)$.

The associated correspondence $\Omega \colon G \to H$ is the space $H = G \ltimes X$. We again consider the bar resolutions $\Z[G^{\bullet + 1} ] \to \Z[G^0]$ and $\Z[H^{\bullet +1}] \to \Z[X]$. For each $n \geq 0$ the induced module $\Ind_\Omega \Z[H^{n+1}]$ has underlying abelian group $\Z[H^{n+1}]$, with the $G$-module structure coming from the action $G \acts H^{n+1}$. The proper $G$-equivariant maps $\tau_n \colon G^{n+1} \to H^{n+1}$ induce a chain map
\begin{equation*}
f := \tau_{\bullet+1}^* \colon \Z[G^{\bullet +1}] \to \Ind_{\Omega} \Z[H^{\bullet+1}]
\end{equation*}
over $\tau^* \colon \Z[G^0] \to \Z[X] = \Ind_{\Omega} \Z[X]$. By Corollary \ref{main corollary}, $\H_*(\Omega)$ is induced by the chain map $( \delta_{\Omega} \otimes \id ) \circ f_G \colon \Z[G^{\bullet+1}]_G \to \Z[H^{\bullet+1}]_H$. The coinvariants $(\Ind_\Omega \Z[H^{n+1}])_G$ is given by $\Z[H^n]$, and $\delta_\Omega \otimes \id \colon \Z[H^n] \to \Z[H^n]$ is simply the identity. The chain map $f_G \colon  (\Z[G^{\bullet +1}])_G \to (\Ind_{\Omega} \Z[H^{\bullet+1}])_G$ is given at $n$ by $\tau_n^*  \colon \Z[G^n] \to \Z[H^n]$. Therefore $(\delta_\Omega \otimes \id ) \circ f_G = \tau_\bullet^*$, and we deduce that $\H_*(\Omega) = \H_*(\tau)$.
\end{example}

\begin{example}[Inverse semigroup homology isomorphism]\label{inverse semigroup homology isomorphism}
Let $S$ be an inverse semigroup, let $E$ be its idempotent semilattice with $E^\times = E \setminus \{ 0 \}$ and let $\hat E$ be the space of filters on $E$. We associate to $S$ both the discrete groupoid $S \ltimes E^\times$ and the universal groupoid $G_S = S \ltimes \hat E$. Then there is a proper correspondence $\Omega_S \colon S \ltimes E^\times \to G_S$ with bispace 
\[ \Omega_S = \bigsqcup_{e \in E^\times} \left\{ [s,\chi] \in S \ltimes \hat E \suchthat s \bcdot \chi \in U_e \right\} \]
such that $\H_*(\Omega_S) \colon \H_*(S \ltimes E^\times) \to \H_*(G_S)$ is an isomorphism. The isotropy of the discrete groupoid $S \ltimes E^\times$ at $e \in E^\times$ is given by the stabiliser subgroup $S_e = \{ s \in S \suchthat s^* s = e = s s^* \}$. The homology of the universal groupoid $G_S$ therefore reduces to a direct sum
\[ \H_*(G_S) \cong \H_*(S \ltimes E^\times) \cong \bigoplus_{[e] \in S \backslash E^\times} \H_*(S_e) \]
of the group homology of the stabiliser subgroups over the orbit space $S \backslash E^\times$.
\begin{proof}
Set $G = S \ltimes E^\times$ and $H = G_S$. Consider the space $Z = \bigsqcup_{e \in E^\times} U_e$ equipped with the diagonal action of $S$ and let $L = S \ltimes Z$. The local homeomorphism $Z \to \hat E$ sending $(e,\chi)$ to $\chi$ induces an \'etale homomorphism $\psi \colon L \to H$. This induces a chain map $\psi^\bullet_* \colon \Z[L^\bullet] \to \Z[H^\bullet]$ of Matui's chain complexes. Furthermore, $L$ is the transformation groupoid of an action $G \acts Z$ with a proper anchor map $\tau \colon Z \to E^\times$ sending $(e,\chi)$ to $e$. This induces a chain map $\tau_\bullet^* \colon \Z[G^\bullet] \to \Z[L^\bullet]$. Setting $\Omega_S \colon G \to H$ to be the composition of the associated correspondences, the map $\H_*(\Omega_S)$ is induced by the chain map $\psi^\bullet_* \circ \tau_\bullet^* \colon \Z[G^\bullet] \to \Z[H^\bullet]$ by Examples \ref{map in homology from etale homomorphism} and \ref{map in homology from action correspondence}. For each $n \geq 0$ and $x \in G^n$, the preimage $\tau_n^{-1}(x)$ is a compact Hausdorff open set on which $\psi^n$ is injective. The image $V_x \defeq \psi^n(\tau_n^{-1}(x))$ is a compact open Hausdorff subset of $H^n$ on which the source map $s_n \colon H^n \to \hat E$ is an isomorphism $V_x \cong U_{s_n(x)}$. Thus for $x \in G^n$ we have
\begin{align*}
\psi^n_* \circ \tau_n^* \colon \Z[G^n] & \to \Z[H^n] \\
 \chi_{\{x\}} & \mapsto \chi_{V_x}.
\end{align*}
We claim that for each $n$ this is an isomorphism. For surjectivity, it suffices to span $\Z[V_x]$ for each $x \in G^n$, as $\{V_y \suchthat y \in G^n \}$ is an open cover of $H^n$. The set $\{ V_y \suchthat y \in G^n, \; V_y \leq V_x \}$ is mapped through $s_n$ to $\{ U_e \suchthat e \in E^\times, \; e \leq s_n(x) \}$. Intersections of the form $\bigcap_{e \in I} U_{e} \cap \bigcap_{e \in J} (U_{s_n(x)} \setminus U_{e})$ for finite subsets $I$ and $J$ of $\{ e \in E^\times \suchthat e \leq s_n(x) \}$ form a basis for the topology of $U_{s_n(x)}$ and therefore $\algspan\{ \chi_{U_e} \suchthat e \leq s_n(x) \} = \Z[U_{s_n(x)}]$ as it is closed under products and contains the indicators of the complements. It follows that $\{\chi_{V_y} \suchthat y \in G^n, \; V_y \subseteq V_x \}$ spans $\Z[V_x]$.

We turn to injectivity. For each $e \in E^\times$ consider the principal filter $e^\uparrow = \{ f \in E \suchthat e \leq f \} \in \hat E$. For each $x \in G^n$ there is a unique $x^\uparrow \in V_x$ with $s_n(x^{\uparrow}) = s_n(x)^\uparrow$. The map $x \mapsto x^\uparrow$ induces a linear map $\Z[  H^n ] \to C_b(G^n, \Z)$. For each $x \in G^n$, the indicator $\chi_{V_x}$ is sent to the indicator $\chi_{x^\downarrow}$ on the set $x^\downarrow = \{ y \in G^n \suchthat V_y \subseteq V_x \}$. Now let $J \subseteq G^n$ be finite and suppose there are integers $(a_x)_{x \in J}$ such that $\sum_{x \in J} a_x \chi_{x^\downarrow} = 0$. Then $a_x = 0$ for any $x \in J$ with $V_x$ maximal, and therefore $a_x = 0$ for each $x \in J$. It follows that $\{ \chi_{V_x} \suchthat x \in G^n \}$ is linearly independent.
\end{proof}
\end{example}

\begin{rmk}
Steinberg constructs an isomorphism between the groupoid ring $\Z[G_S]$ and the inverse semigroup ring $\Z[S]$ \cite[Theorem 6.3]{Steinberg10}, from which the above proof that $\H_*(G_S) \cong \H_*(S \ltimes E^\times)$ takes heavy inspiration. The feature of Example \ref{inverse semigroup homology isomorphism} we would like to highlight is how this isomorphism in homology is induced by an explicit étale correspondence.
\end{rmk}


\begin{thebibliography}{BDGW23}

\bibitem[Alb15]{Albandik15}
Suliman Albandik.
\newblock {\em A colimit construction for groupoids}.
\newblock 2015.
\newblock Thesis (Ph.D.)--{G}eorg-{A}ugust-{U}niversit\"at {G}\"ottingen.

\bibitem[AKM22]{AKM22}
Celso Antunes, Joanna Ko, and Ralf Meyer.
\newblock The bicategory of groupoid correspondences.
\newblock {\em New York J. Math.}, 28:1329--1364, 2022.


\bibitem[BDGW23]{BDGW23}
Christian B\"onicke, Cl\'ement Dell'Aiera, James Gabe, and Rufus Willett.
\newblock Dynamic asymptotic dimension and {M}atui's {HK} conjecture.
\newblock {\em Proc. Lond. Math. Soc. (3)}, 126(4):1182--1253, 2023.


\bibitem[CM00]{CraMoe00}
Marius Crainic and Ieke Moerdijk.
\newblock A homology theory for \'{e}tale groupoids.
\newblock {\em J. Reine Angew. Math.}, 521:25--46, 2000.

\bibitem[Dee23]{Deeley21}
Robin~J. Deeley.
\newblock A counterexample to the {HK}-conjecture that is principal.
\newblock {\em Ergodic Theory Dynam. Systems}, 43(6):1829--1846, 2023.

\bibitem[FKPS19]{FKPS19}
Carla Farsi, Alex Kumjian, David Pask, and Aidan Sims.
\newblock Ample groupoids: equivalence, homology, and {M}atui's {HK}
  conjecture.
\newblock {\em M\"{u}nster J. Math.}, 12(2):411--451, 2019.

\bibitem[Hol14]{Holkar14}
Rohit~Dilip Holkar.
\newblock Topological construction of {$\mathrm{C}^\ast$}-correspondences for
  groupoid {$\mathrm{C}^\ast$}-algebras.
\newblock 2014.
\newblock Ph.D. Thesis, Georg-August-Universit\"at G\"ottingen.

\bibitem[Hol17a]{Holkar17b}
Rohit~Dilip Holkar.
\newblock Composition of topological correspondences.
\newblock {\em J. Operator Theory}, 78(1):89--117, 2017.

\bibitem[Hol17b]{Holkar17a}
Rohit~Dilip Holkar.
\newblock Topological construction of {$\mathrm{C}^\ast$}-correspondences for
  groupoid {$\mathrm{C}^\ast$}-algebras.
\newblock {\em J. Operator Theory}, 77(1):217--241, 2017.

\bibitem[KS06]{KasSch06}
Masaki Kashiwara and Pierre Schapira.
\newblock {\em Categories and sheaves}, volume 332 of {\em Grundlehren der
  mathematischen Wissenschaften [Fundamental Principles of Mathematical
  Sciences]}.
\newblock Springer-Verlag, Berlin, 2006.

\bibitem[Li22]{Li22}
Xin Li.
\newblock Ample groupoids, topological full groups, algebraic {K}-theory
  spectra and infinite loop spaces.
\newblock 2022.
\newblock Preprint, arXiv:2209.08087v1.

\bibitem[Mat12]{Matui12}
Hiroki Matui.
\newblock Homology and topological full groups of \'{e}tale groupoids on
  totally disconnected spaces.
\newblock {\em Proc. Lond. Math. Soc. (3)}, 104(1):27--56, 2012.

\bibitem[Mat15]{Matui15}
Hiroki Matui.
\newblock Topological full groups of one-sided shifts of finite type.
\newblock {\em J. Reine Angew. Math.}, 705:35--84, 2015.

\bibitem[Mat16]{Matui16}
Hiroki Matui.
\newblock \'{E}tale groupoids arising from products of shifts of finite type.
\newblock {\em Adv. Math.}, 303:502--548, 2016.

\bibitem[Mil24]{Miller23a}
Alistair Miller.
\newblock Functors between equivariant {K}asparov categories from \'{e}tale
  groupoid correspondences.
\newblock {\em J. Funct. Anal.}, 287(11):Paper No. 110623, 2024.

\bibitem[NO21a]{NylOrt21a}
Petter Nyland and Eduard Ortega.
\newblock Katsura-{E}xel-{P}ardo groupoids and the {AH} conjecture.
\newblock {\em J. Lond. Math. Soc. (2)}, 104(5):2240--2259, 2021.

\bibitem[NO21b]{NylOrt21b}
Petter Nyland and Eduard Ortega.
\newblock Matui's {AH} conjecture for graph groupoids.
\newblock {\em Doc. Math.}, 26:1679--1727, 2021.

\bibitem[Ort20]{Ortega20}
Eduard Ortega.
\newblock The homology of the {K}atsura-{E}xel-{P}ardo groupoid.
\newblock {\em J. Noncommut. Geom.}, 14(3):913--935, 2020.

\bibitem[PY22]{ProYam22}
Valerio Proietti and Makoto Yamashita.
\newblock Homology and {K}-theory of dynamical systems {I}. {T}orsion-free
  ample groupoids.
\newblock {\em Ergodic Theory Dynam. Systems}, 42(8):2630--2660, 2022.

\bibitem[Rot09]{Rotman09}
Joseph~J. Rotman.
\newblock {\em An introduction to homological algebra}.
\newblock Universitext. Springer, New York, second edition, 2009.

\bibitem[Sca20]{Scarparo20}
Eduardo Scarparo.
\newblock Homology of odometers.
\newblock {\em Ergodic Theory Dynam. Systems}, 40(9):2541--2551, 2020.

\bibitem[Sca23]{Scarparo22}
Eduardo Scarparo.
\newblock {AH} conjecture for {C}antor minimal dihedral systems.
\newblock {\em Math. Scand.}, 129(2):284--298, 2023.

\bibitem[Ste10]{Steinberg10}
Benjamin Steinberg.
\newblock A groupoid approach to discrete inverse semigroup algebras.
\newblock {\em Adv. Math.}, 223(2):689--727, 2010.

\bibitem[Ste14]{Steinberg14}
Benjamin Steinberg.
\newblock Modules over \'{e}tale groupoid algebras as sheaves.
\newblock {\em J. Aust. Math. Soc.}, 97(3):418--429, 2014.

\bibitem[Wei94]{Weibel94}
Charles~A. Weibel.
\newblock {\em An introduction to homological algebra}, volume~38 of {\em
  Cambridge Studies in Advanced Mathematics}.
\newblock Cambridge University Press, Cambridge, 1994.

\end{thebibliography}
\end{document}